\documentclass[11pt,epsfig]{article}
\usepackage{graphicx,wrapfig}
\usepackage{flafter}
\usepackage{leftidx}
\usepackage{mathrsfs}
\usepackage{amssymb,amsmath}
\usepackage{bbding}
\usepackage{fancyhdr}
\usepackage{mathrsfs}
\usepackage{color}
\usepackage{stmaryrd}
\usepackage{amsfonts}
\usepackage{latexsym}
\usepackage{psfrag}
\usepackage{graphicx,subfigure}
\usepackage{fancyhdr,graphicx}
\usepackage{multicol}
\usepackage{dsfont}
\usepackage{bbm}
\usepackage{booktabs}
\usepackage[center]{caption2}
\usepackage{cite}
\usepackage{epstopdf}
\usepackage{diagbox}
\usepackage{booktabs}
\usepackage [latin1]{inputenc}
\usepackage{multirow}
\usepackage{enumerate}
\usepackage{enumitem}
\usepackage{algpseudocode,algorithm,algorithmicx}

\allowdisplaybreaks

\date{}
\newtheorem{theorem}{Theorem}[section]
\newtheorem{lemma}{Lemma}[section]

\numberwithin{equation}{section}

\newcommand{\rmk}{\mathtt{k}}
\newcommand{\zd}{\,\mathrm{d}}

\newcommand{\defeq}{:=}
\newcommand{\diff}{\triangledown_{\tau}}
\newcommand{\abs}[1]{\left|#1\right|}
\newcommand{\absb}[1]{\big|#1\big|}

\newcommand{\bra}[1]{\left(#1\right)}
\newcommand{\brab}[1]{\big(#1\big)}
\newcommand{\braB}[1]{\Big(#1\Big)}
\newcommand{\brat}[1]{(#1)}

\newcommand{\kbrab}[1]{\big[#1\big]}

\newcommand{\myinner}[1]{\left\langle#1\right\rangle}
\newcommand{\myinnerb}[1]{\big\langle#1\big\rangle}

\newcommand{\mynorm}[1]{\left\|#1\right\|}
\newcommand{\mynormb}[1]{\big\|#1\big\|}

\newcommand{\myvec}[1]{\boldsymbol{#1}}

\title{A new discrete energy technique for multi-step backward difference formulas}
\author{Hong-lin Liao\thanks{ORCID 0000-0003-0777-6832;
Department of Mathematics,
Nanjing University of Aeronautics and Astronautics,
Nanjing 211106, P. R. China. E-mails: liaohl@nuaa.edu.cn and liaohl@csrc.ac.cn.
This author's work is supported by NSF of China under grant number 12071216.}
\qquad Tao Tang\thanks{Division of Science and Technology, BNU-HKBU United International College,
    Zhuhai, Guangdong, China, and International Center for Mathematics, Southern
    University of Science and Technology, Shenzhen, Guangdong, China.
    Email: ttang@uic.edu.cn. This author's work is partially supported by NSF of China (under grant numbers 11731006 and K20911001) and the NNW2018-ZT4A06 project.}
\qquad Tao Zhou\thanks{LSEC, Institute of Computational Mathematics and Scientific/Engineering Computing,
Academy of Mathematics and Systems Science, Chinese Academy of Sciences, Beijing, 100190, China. Email: tzhou@lsec.cc.ac.cn. The work of T. Zhou is supported by NSF of China (under grant numbers 11822111, 11688101), the Strategic Priority Research Program of Chinese Academy of Sciences (Grant No. XDA25000404), the Science Challenge Project (No. TZ2018001), and Youth Innovation Promotion Association of CAS.}
}
\date{\today}

\begin{document}

\maketitle

\begin{abstract}
The backward differentiation formula (BDF) is a useful family of implicit methods for the numerical integration of stiff differential equations. It is well noticed that the stability and convergence of the $A$-stable BDF1 and BDF2 schemes for parabolic equations can be directly established by using the standard discrete energy analysis. However, such classical analysis technique seems not directly applicable to the BDF-$\rmk$ schemes for $3\le \rmk\le 5$. To overcome the difficulty, a powerful analysis tool based on the Nevanlinna-Odeh multiplier technique [Numer. Funct. Anal. Optim., 3:377-423, 1981] was developed by Lubich et al. [IMA J. Numer. Anal., 33:1365-1385, 2013]. In this work, by using the so-called discrete orthogonal convolution kernels technique, we will recover the classical energy analysis so that the stability and convergence of the BDF-$\rmk$ schemes for $3\le \rmk\le 5$ can be established. One of the theoretical advantages of our analysis technique is that less spacial regularity requirement is needed on the initial data.
\\[1ex]
\end{abstract}

\section{Introduction}
\setcounter{equation}{0}

The backward differentiation formula (BDF) is a family of implicit methods for the numerical integration of stiff differential equations. They are linear multistep methods that approximate the derivative of the unknown function using information from already computed time points, thereby increasing the accuracy of the approximation. These methods are especially used for the solution of stiff differential equations whose numerical stability is indicated by the region of absolute stability. More precisely, if the region of stability contains the left half of the complex plane, then the method is said to be $A$-stable. It is known that backward differentiation methods with an order greater than 2 cannot be $A$-stable, i.e., only the first-order and second-order backward differentiation formulas (i.e., BDF1 and BDF2) are $A$-stable. For parabolic equations, it is also well-known that the energy stability and convergence of $A$-stable BDF1 and BDF2 methods can be established by using the standard discrete energy analysis, see, e.g., \cite{Thomee:2006}.

However, this standard analysis technique is not directly applicable to higher order  BDF schemes. This results in several remedies to recover  the $L^2$-norm
stability and convergence for the non-$A$-stable $\rmk$-step backward difference formulas with $3\le \rmk\le5$. It is particularly noted that due to the seminal work of Lubich et al. \cite{LubichMansourVenkataraman:2013}, the Nevanlinna-Odeh multiplier technique \cite{NevanlinnaOdeh:1981} has been successfully used for this purpose, see e.g. \cite{Akrivis:2015,AkrivisKatsoprinakis:2016,AkrivisChenYuZhou:2021,AkrivisLubich:2015} and references therein. The key idea of such a multiplier technique relies on the equivalence between $A$-stability and $G$-stability
 of Dahlquist \cite{Dahlquist:1978}. Another useful tool for the numerical analysis of  BDF-$\rmk$ schemes is the telescope formulas by Liu \cite{Liu:2013}, which is also based on the Dahlquist's $G$-stability theory \cite{Dahlquist:1978}.

We have a natural question: is there a straightforward discrete energy analysis for the  BDF-$\rmk$ schemes with $3\le \rmk\le5$?
The aim of this work is to provide a definite answer by introducing a novel yet straightforward discrete energy method based on the so-called discrete orthogonal convolution (DOC) kernel technique \cite{LiaoZhang:2020}.
To this end, we consider the linear reaction-diffusion problem in a bounded convex domain $\Omega$,
\begin{align}\label{eq: diffusion problem}
\partial_tu-\varepsilon\Delta u&=\beta (x,t)u+f(t,x), \quad x\in\Omega, \quad 0<t<T,
\end{align}
subject to the Dirichlet boundary condition~$u=0$ on a smooth boundary~$\partial\Omega.$ The initial data is set to be $u(0,x)=u_0(x).$
We assume that the diffusive coefficient $\varepsilon>0$ is a constant
and the reaction coefficient $\beta (x,t)$ is smooth  and bounded by $\beta ^{*}>0$.

Let $t_k=k\tau$ be an uniform discrete time-step with $\tau\defeq T/N$. For any discrete time sequence $\{v^n\}_{n=0}^N$,
we denote
\[
\diff v^n:=v^n-v^{n-1}, \quad \partial_{\tau}v^n:=\diff v^n/\tau.
\]
 For a fixed index $3\le \rmk\le 5$, we shall view the BDF-$\rmk$  formula $D_{\rmk}v^n$ as a discrete convolution summation as follows
\begin{align}\label{def: BDF2-Formula}
D_{\rmk}v^n:=\frac1{\tau}\sum_{k=1}^nb_{n-k}^{(\rmk)}\diff v^k,
\quad  n\ge \rmk,
\end{align}
where the associated BDF-$\rmk$ kernels $b_{j}^{(\rmk)}$ (vanish if $j\ge \rmk$, see Table \ref{table: BDF-k kernels}) are generated by
\begin{align}\label{def: BDF-k generating function}
\sum_{\ell=1}^{\rmk}\frac{1}{\ell}(1-\zeta)^{\ell-1}
=\sum_{\ell=0}^{\rmk-1}b_{\ell}^{(\rmk)}\zeta^\ell, \quad 3\le \rmk\le 5.
\end{align}

\begin{table}[htb!]
\begin{center}\label{table: BDF-k kernels}
\caption{The BDF-$\rmk$ kernels $b_{j}^{(\rmk)}$ generated by \eqref{def: BDF-k generating function}}
\vspace*{0.3pt}
\def\temptablewidth{0.7\textwidth}
{\rule{\temptablewidth}{0.5pt}}
\begin{tabular*}{\temptablewidth}{@{\extracolsep{\fill}}cccccc}
   BDF-$\rmk$      &$b_{0}^{(\rmk)}$     &$b_{1}^{(\rmk)}$     &$b_{2}^{(\rmk)}$
   &$b_{3}^{(\rmk)}$  &$b_{4}^{(\rmk)}$ \\  \midrule
  $\rmk=2$   &$\frac{3}{2}$   &$-\frac{1}{2}$   &&&\\[3pt]
  $\rmk=3$   &$\frac{11}{6}$  &$-\frac{7}{6}$   &$\frac{1}{3}$   &&\\[3pt]
  $\rmk=4$   &$\frac{25}{12}$ &$-\frac{23}{12}$ &$\frac{13}{12}$ &$-\frac{1}{4}$ &\\[3pt]
  $\rmk=5$   &$\frac{137}{60}$&$-\frac{163}{60}$&$\frac{137}{60}$&$-\frac{21}{20}$ &$\frac{1}{5}$
\end{tabular*}
{\rule{\temptablewidth}{0.5pt}}
\end{center}
\end{table}	
	
To make our idea clear, the initial data $ u^1$, $ u^2$, $\cdots$, $ u^{\rmk-1}$ for
the multi-step BDF-$\rmk$ schemes ($3\le \rmk\le 5$) are
assumed to be available. Without loss of generality, we consider the time-discrete solution, $u^k(x)\approx u(t_k,x)$
for~$x\in\Omega$, defined by the following implicit multi-step BDF scheme
\begin{align}\label{eq: time-discrete IBVP}
D_{\rmk}u^{j}
	=\varepsilon\Delta u^{j}
        +\beta ^j u^{j}+f^j,\quad  \rmk\le j\le N,
\end{align}
where $f^j(x)=f(t_j,x)$. The weak form of \eqref{eq: time-discrete IBVP} reads
\begin{align}\label{eq: weak time-discrete IBVP}
\myinnerb{D_{\rmk}u^{j},w}
	+\varepsilon\myinnerb{\nabla u^{j},\nabla w}
        =\myinnerb{\beta ^j u^{j},w}+\myinnerb{f^j,w}\quad\text{for $\forall\;w\in  H^1_0(\Omega)$ and $\rmk\le j\le N.$}
\end{align}
Here $\myinner{\cdot,\cdot}$ and $\mynorm{\cdot}$ are the $L^2$ inner product and the $L^2$-norm, respectively.

Our new energy analysis for BDF-$\rmk$ schemes with $3\le \rmk\le5$ relies on the discrete orthogonal convolution (DOC) kernels technique developed in \cite{LiaoZhang:2020},
where the BDF2 scheme (with variable steps) was investigated. More precisely, our analysis will be based on an equivalent convolution form of \eqref{eq: weak time-discrete IBVP} using the DOC kernels technique. Below we will derive the equivalent convolution form.
For the discrete BDF-$\rmk$ kernels $b_{j}^{(\rmk)}$ generated by \eqref{def: BDF-k generating function},
the corresponding DOC-$\rmk$ kernels $\theta_{j}^{(\rmk)}$ are defined  recursively as:
\begin{align}\label{def: DOC-Kernels}
\theta_{0}^{(\rmk)}:=\frac{1}{b_{0}^{(\rmk)}}
\quad \mathrm{and} \quad
\theta_{n-j}^{(\rmk)}:=-\frac{1}{b_{0}^{(\rmk)}}
\sum_{\ell=j+1}^n\theta_{n-\ell}^{(\rmk)}b_{\ell-j}^{(\rmk)}\quad \text{for $j=n-1,n-2,\cdots,\rmk$.}
\end{align}
Here and hereafter, we set $\sum_{k=i}^{j}\cdot=0$ whenever $i>j$.
It is easy to check that the following \emph{discrete orthogonal convolution identity} holds \cite{LiaoZhang:2020,LiaoTangZhou:2020doc}:
\begin{align}\label{eq: orthogonal identity}
\sum_{\ell=j}^{n}\theta_{n-\ell}^{(\rmk)}b^{(\rmk)}_{\ell-k}\equiv\delta_{nj}
\quad\text{for $\forall\;\rmk\leq j\le n$,}
\end{align}
where $\delta_{nk}$ is the Kronecker delta. Thus, by exchanging the summation index, one gets
\begin{align*}
\sum_{j=\rmk}^{n}\theta_{n-j}^{(\rmk)}
\sum_{\ell=\rmk}^{j}b_{j-\ell}^{(\rmk)}\diff u^\ell
=&\,\sum_{\ell=\rmk}^{n}\diff u^\ell\sum_{j=\ell}^{n}\theta_{n-j}^{(\rmk)}b_{j-\ell}^{(\rmk)}
=\diff u^n, \quad \rmk\le n\le N.
\end{align*}
Then by acting the associated DOC kernels $\theta_{n-j}^{(\rmk)}$ on the BDF-k formula $D_{\rmk}$, we obtain
\begin{align}\label{Dis: DOC action BDF formula Dk}
\sum_{j=\rmk}^{n}\theta_{n-j}^{(\rmk)}D_{\rmk} u^j
=&\,\frac1{\tau}\sum_{j=\rmk}^{n}\theta_{n-j}^{(\rmk)}
\sum_{\ell=1}^{\rmk-1}b_{j-\ell}^{(\rmk)}\diff u^\ell
+\frac1{\tau}\sum_{j=\rmk}^{n}\theta_{n-j}^{(\rmk)}
\sum_{\ell=\rmk}^{j}b_{j-\ell}^{(\rmk)}\diff u^\ell\nonumber\\
\triangleq&\, \frac1{\tau}u_{\mathrm{I}}^{(\rmk,n)}+\partial_{\tau} u^n\qquad\text{for $\rmk\le n\le N$,}
\end{align}
where $ u_{\mathrm{I}}^{(\rmk,n)}$ represents the starting effects
on the numerical solution at $t_n,$ i.e.,
\begin{align}\label{Dis: initial effect -BDF formula}
 u_{\mathrm{I}}^{(\rmk,n)}:=\sum_{\ell=1}^{\rmk-1}\diff u^\ell
\sum_{j=\rmk}^{n}\theta_{n-j}^{(\rmk)}b_{j-\ell}^{(\rmk)}\qquad\text{for $n\ge\rmk$.}
\end{align}
Now, by acting the associated DOC-k kernels $\theta_{n-j}^{(\rmk)}$ on the time-discrete problem
\eqref{eq: weak time-discrete IBVP}, we use \eqref{Dis: DOC action BDF formula Dk}
and \eqref{Dis: initial effect -BDF formula} to obtain
\begin{align}\label{eq: DOC weak time-discrete IBVP}
\myinnerb{\partial_{\tau} u^{j},w}+\varepsilon\sum_{\ell=\rmk}^j\theta_{j-\ell}^{(\rmk)}\myinnerb{\nabla u^{\ell},\nabla w}
        =&\,-{1\over \tau} \myinnerb{u_{\mathrm{I}}^{(\rmk,j)},w}\,+\sum_{\ell=\rmk}^j\theta_{j-\ell}^{(\rmk)}\myinnerb{\beta ^{\ell}u^{\ell},w}
        +\sum_{\ell=\rmk}^j\theta_{j-\ell}^{(\rmk)}\myinnerb{f^\ell,w}\nonumber\\
        &\,\qquad\text{for $\forall\;w\in  H^1_0(\Omega)$ and $\rmk\le j\le N.$}
\end{align}
This convolution formulation will be the starting point of our energy technique, and will lead
to much more concise $L^2$-norm estimates than those in previous works.

Note that, the DOC-$\rmk$ kernels define a \textit{reversible discrete transform}
between the original form \eqref{eq: weak time-discrete IBVP}
and the convolution form \eqref{eq: DOC weak time-discrete IBVP}.
By following the proof of \cite[Lemma 2.1]{LiaoTangZhou:2020doc}, it is easy to check
 \begin{align}\label{eq: mutual orthogonal identity}
 \sum_{\ell=j}^nb_{n-\ell}^{(\rmk)}\theta_{\ell-j}^{(\rmk)}\equiv \delta_{mk}\quad\text{for $\rmk\le j\le n$.}
   \end{align}

Similar as the classical discrete $L^2$-norm analysis, we can now take $w=2\tau u^j$ in \eqref{eq: DOC weak time-discrete IBVP}
and sum up from $j=\rmk$ to $n$ to obtain
\begin{align}\label{eq: DOC IBVP energy form}
\mynormb{u^{n}}^2-\mynormb{u^{\rmk-1}}^2
        \le&\,-2\sum_{j=\rmk}^n\myinnerb{ u_{\mathrm{I}}^{(\rmk,j)},u^j} -2\varepsilon\tau\sum_{j=\rmk}^n\sum_{\ell=\rmk}^j\theta_{j-\ell}^{(\rmk)}\myinnerb{\nabla u^{\ell},\nabla u^j}\nonumber\\
        &\,        +2\tau\sum_{j=\rmk}^n\sum_{\ell=\rmk}^j\theta_{j-\ell}^{(\rmk)}\myinnerb{\beta ^{\ell}u^{\ell},u^j}
+2\tau\sum_{j=\rmk}^n\sum_{\ell=\rmk}^j\theta_{j-\ell}^{(\rmk)}\myinnerb{f^\ell,u^j},
\end{align}
where the term $\sum_{j=\rmk}^n \|u^j-u^{j-1}\|^2 $ has been dropped in the left hand side. Consequently, we need to carefully handle the right hand side of \eqref{eq: DOC IBVP energy form}, which consists the following issues:
\begin{itemize}
  \item Positive definiteness of the DOC kernels $\theta_{j}^{(\rmk)}$ (see Section 2.1);
  \item Decay estimates of the DOC kernels $\theta_{j}^{(\rmk)}$ (see Section 2.2);
  \item Decay estimates of the initial term $u_{\mathrm{I}}^{(\rmk,j)}$ (see Section 2.3).
\end{itemize}

By doing this, we can finally present concise stability and error estimates of \eqref{eq: time-discrete IBVP}
for the linear reaction-diffusion \eqref{eq: diffusion problem}. More precisely, we show in
Theorem \ref{thm: general L2 norm stability} that,
if the time-step size $\tau\le(7-\rmk)/(7\rho_{\rmk}\beta ^{*})$,
the time-discrete solution $u^n$ is unconditionally stable in $L^2$-norm:
\begin{align*}
\mynormb{u^{n}}\le&\, \frac{7\rho_{\rmk}}{7-\rmk}\exp\braB{\frac{7\rho_{\rmk}}{7-\rmk}\beta^{*}t_{n-\rmk}}
\braB{c_{\mathrm{I},\rmk}\sum_{\ell=0}^{\rmk-1}\mynormb{u^\ell}
+\sum_{\ell=\rmk}^n\tau\mynormb{f^\ell}}\quad\text{for $\rmk\le n\le N$},
\end{align*}
where the constants $\rho_{\rmk}$ and $c_{\mathrm{I},\rmk}$ are defined in Lemmas \ref{lem: decaying estimates DOC}
and \ref{lem: Estimates of initial terms uI}, respectively.
This is followed by Theorem \ref{thm: general L2 norm converegnce} which presents
a concise $L^2$-norm error estimate for the BDF-$\rmk$ schemes:
\begin{align*}
\mynormb{u(t_n)-u^{n}}\le&\,
\frac{7\rho_{\rmk}c_{\mathrm{I},\rmk}}{7-\rmk}
\exp\braB{\frac{7\rho_{\rmk}\beta^{*}t_{n-\rmk}}{7-\rmk}}
\bigg(\sum_{\ell=0}^{\rmk-1}\mynormb{u(t_{\ell})-u^{\ell}}
+C_ut_{n-\rmk}\tau^{\rmk}\bigg).
\end{align*}

The paper is organized as follows. Section 2 contains several properties for the DOC-$\rmk$ kernels which will be useful for the stability and convergence analysis. The main results outlined above will be proved in Section 3. Some concluding remarks will be given in the last section.

\section{Preliminary results}
In this section, we present several preliminary results which will be used for proving our main results in Section 3.

\subsection{Positive definiteness of DOC-$\rmk$ kernels}
\label{section: DOC positive definiteness}
\setcounter{equation}{0}

By using the mutual orthogonal identities \eqref{eq: orthogonal identity} and \eqref{eq: mutual orthogonal identity},
we have the following result on the positive definiteness, see \cite[Lemma 2.1]{LiaoTangZhou:2020doc}.
\begin{lemma}\label{lem: equvielant positive definite DOC}
The discrete kernels $b^{(\rmk)}_{j}$ in \eqref{def: BDF-k generating function} are positive (semi-)definite
if and only if the associated DOC kernels $\theta_{j}^{(\rmk)}$ in \eqref{def: DOC-Kernels}  are positive (semi-)definite.
\end{lemma}

It remains to study the positive definiteness of discrete BDF-$\rmk$ kernels. To this end, we
introduce the Toeplitz form $\mathrm{T}_m=(\mathrm{t}_{ij})_{m\times m}$, where the entries
$\mathrm{t}_{ij}=\mathrm{t}_{i-j}$ ($i,j=1,2,\cdots, m$) are constants along the diagonal of $\mathrm{T}_m$.
Let $\mathrm{t}_{k}$ be the Fourier coefficients of the trigonometric polynomial $\mathrm{f}$, i.e.,
$$\mathrm{t}_{k}=\frac{1}{2\pi}\int_{-\pi}^{\pi}\mathrm{f}(x)e^{-\imath k x}\zd x, \quad k=1-m,2-m, \cdots, m-1,$$
where $\imath:=\sqrt{-1}$ is the complex unit. Then
\begin{align}\label{eq: generating function Toeplitz}
\mathrm{f}(x)=\sum_{k=1-m}^{m-1}\mathrm{t}_{k}e^{\imath k x}
\quad \text{is called the generating function of $\mathrm{T}_m$}.
\end{align}
The following Grenander-Szeg\"{o} theorem \cite[pp. 64--65]{GrenanderSzego:2001} shows the relationship
between the eigenvalues of $\mathrm{T}_m$ and the generating function $\mathrm{f}$.

\begin{lemma}\label{lem: Grenander Szego theorem}
Let $\mathrm{T}_m=(\mathrm{t}_{i-j})_{m\times m}$ be the Toeplitz matrix generated by the function
$\mathrm{f}$ defined in \eqref{eq: generating function Toeplitz}.
Then the smallest eigenvalue $\lambda_{\min}\brab{\mathrm{T}_m}$ of $\mathrm{T}_m$
and the largest one $\lambda_{\max}\brab{\mathrm{T}_m}$ can be bounded by
$$\mathrm{f}_{\min}\le\lambda_{\min}\brab{\mathrm{T}_m}\le\lambda_{\max}\brab{\mathrm{T}_m}\le \mathrm{f}_{\max},$$
where $\mathrm{f}_{\min}$ and $\mathrm{f}_{\max}$ denote the minimum and maximum values of $\mathrm{f}(x)$, respectively.
Specially, the Toeplitz matrix $\mathrm{T}_m$ is positive definite if $\mathrm{f}_{\min}>0$.
\end{lemma}

Now, for the BDF-$\rmk$ formula, we consider the following matrices of order $m:=n-\rmk+1:$
\begin{align}\label{matrix: Bk BDF-k formula}
B_{\rmk,l}:=\left(
 \begin{array}{cccccc}
    b_0^{(\rmk)} & && &&\\
     b_1^{(\rmk)} &b_0^{(\rmk)} & &&&\\
      \vdots            &    \ddots          & \ddots & &    &\\
      b_{\rmk-1}^{(\rmk)}&\cdots &b_1^{(\rmk)}&b_0^{(\rmk)}&&\\
      &\ddots&\cdots&b_1^{(\rmk)}&b_0^{(\rmk)}&\\
      &&b_{\rmk-1}^{(\rmk)}&\cdots&b_1^{(\rmk)}&b_0^{(\rmk)}\\
  \end{array}
\right)_{m\times m}\quad\text{and}\quad B_{\rmk}:=B_{\rmk,l}+B_{\rmk,l}^T,
\end{align}
where $3\le \rmk \le 5$ and the index $n\ge\rmk$.
According to definition \eqref{eq: generating function Toeplitz},
we define the generating function of $B_{\rmk}$ as follows,
\begin{align}\label{eq: generating function Bk}
\mathrm{g}^{(\rmk)}(\varphi)=2\sum_{j=0}^{\rmk-1}b_j^{(\rmk)}\cos(j\varphi).
\end{align}
Consequently, Lemma \ref{lem: Grenander Szego theorem} implies the following result.

\begin{lemma}\label{lem: Max-min Condition positive definite}
Let $B_{\rmk}$ be the real symmetric matrix generated by the function
$\mathrm{g}^{(\rmk)}$ defined in \eqref{eq: generating function Bk}.
Then the smallest eigenvalue $\lambda_{\min}\brab{B_{\rmk}}$
and the largest one $\lambda_{\max}\brab{B_{\rmk}}$ can be bounded by
$$\mathrm{g}^{(\rmk)}_{\min}\le\lambda_{\min}\brab{B_{\rmk}}\le\lambda_{\max}\brab{B_{\rmk}}\le \mathrm{g}^{(\rmk)}_{\max},$$
where $\mathrm{g}^{(\rmk)}_{\min}$ and $\mathrm{g}^{(\rmk)}_{\max}$ denote the minimum and maximum values of $\mathrm{g}^{(\rmk)}(\varphi)$, respectively.
Specially, the real symmetric  matrix $B_{\rmk}$ is positive definite if $\mathrm{g}^{(\rmk)}_{\min}>0$.
\end{lemma}

Now we apply Lemma \ref{lem: Max-min Condition positive definite} to establish
the positive definiteness of the discrete BDF-$\rmk$ kernels $b_j^{(\rmk)}$
for $3\le \rmk \le 5$, by investigating the associated generating functions $\mathrm{g}^{(\rmk)}(\varphi)$.

\begin{lemma}\label{lem: BDF-k minimum eigenvalue}
For the discrete BDF-$\rmk$ kernels $b_{j}^{(\rmk)}$ $(3\le \rmk \le 5)$ defined in \eqref{def: BDF-k generating function}
and any real sequence $\{w_k\}_{k=1}^n$ with n entries, it holds that
\[
2\sum_{m=\rmk}^n w_m \sum_{j=\rmk}^m b_{m-j}^{(\rmk)}w_j
\ge \sigma_\rmk\sum_{k=\rmk}^nw_k^2\quad\text{for $n\ge \rmk$},
\]
where $\sigma_3=95/48\approx1.979$, $\sigma_4=1.628$ and $\sigma_5=0.4776$.
Thus the discrete BDF-$\rmk$ kernels $b_{j}^{(\rmk)}$ for $3\le \rmk \le 5$ are positive definite.
\end{lemma}
\begin{proof} Consider the real symmetric matrix $B_{\rmk}$ in \eqref{matrix: Bk BDF-k formula} of order $m:=n-\rmk+1$. By setting $\myvec{w}:=(w_{\rmk},w_{\rmk+1},\cdots,w_n)^T$
one obtains
\[
2\sum_{m=\rmk}^n w_m \sum_{j=\rmk}^m b_{m-j}^{(\rmk)}w_j
=\myvec{w}^TB_{\rmk}\myvec{w}\ge \lambda_{\min}\brab{B_{\rmk}}\myvec{w}^T\myvec{w} \quad\text{for $n\ge \rmk$}.
\]
Thanks to Lemma \ref{lem: Max-min Condition positive definite},
it remains to prove $\mathrm{g}^{(\rmk)}_{\min}\ge\sigma_{\rmk}.$  The associated generating functions (see Figure \ref{fig: curves generating functions} for the function curves) are listed below:
\begin{itemize}
  \item $\mathrm{g}^{(3)}(\varphi)=\frac{1}{3}\brab{11-7\cos\varphi+2\cos2\varphi},$
  \item $\mathrm{g}^{(4)}(\varphi)=\frac{1}{6}\brab{25-23\cos\varphi+13\cos2\varphi-3\cos3\varphi},$
  \item $\mathrm{g}^{(5)}(\varphi)=\frac{1}{30}\brab{137-163\cos\varphi+137\cos2\varphi-63\cos3\varphi+12\cos4\varphi}.$
\end{itemize}

\begin{figure}[htb!]
\centering
\includegraphics[width=3.3in,height=2.2in]{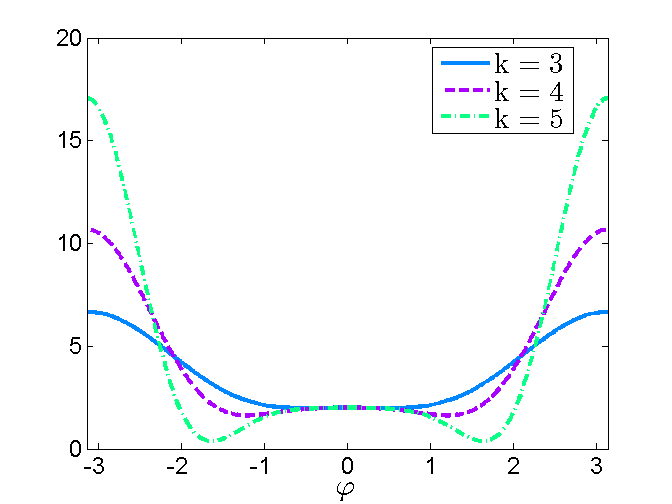}
\caption{Curves of generating functions $\mathrm{g}^{(\rmk)}(\varphi)$ over $\varphi\in[-\pi,\pi]$ for $3\le \rmk \le 5$.}
\label{fig: curves generating functions}
\end{figure}

\textbf{(I) The case $\rmk=3.$}\quad By the formula $\cos2\varphi=2\cos^2\varphi-1$, we get
 \begin{align*}
  \mathrm{g}^{(3)}(\varphi)=\frac{1}{3}\brab{9-7\cos\varphi+4\cos^2\varphi}
  =\frac{4}{3}(\cos\varphi-7/8)^2+\frac{95}{48}.
  \end{align*}
  As desired, $\mathrm{g}^{(3)}_{\min}=\sigma_{3}=95/48\approx 1.97919$.

\textbf{(II) The case $\rmk=4.$}\quad By the formula $\cos3\varphi=4\cos^3\varphi-3\cos\varphi$, we get
\begin{align*}
  \mathrm{g}^{(4)}(\varphi)=\frac{1}{6}\brab{12-14\cos\varphi+26\cos^2\varphi-12\cos^3\varphi}.
  \end{align*}
  Consider a function $Z_4(x)=12-14x+26x^2-12x^3$. The first derivative $Z_4'=-14+52x-36x^2$
  has a unique zero-point $x_{*}=(13-\sqrt{43})/18$ for $x\in[-1,1]$.
  Then $$(Z_4)_{\min}=\min\{Z_4(-1),Z_4(x_{*}),Z_4(1)\}=Z_4(x_{*})=(2656-43\sqrt{43})/243.$$
  Thus, we have
  $$\mathrm{g}^{(4)}(\varphi)\ge \frac1{6}Z_4(x_{*})=\frac{2656-43\sqrt{43}}{1458}\approx1.62828.$$

\textbf{(III) The case $\rmk=5.$}\quad By the formula $\cos4\varphi=8\cos^4\varphi-8\cos^2\varphi+1$, we get
\begin{align*}
  \mathrm{g}^{(5)}(\varphi)=\frac{1}{30}\brab{12+26\cos\varphi+178\cos^2\varphi-252\cos^3\varphi+96\cos^4\varphi}.
  \end{align*}
  Consider the following function $Z_5(x)=12+26x+178x^2-252x^3+96x^4$. The first derivative $Z_5'=26 + 356x - 756 x^2 + 384 x^3$
  has a unique real zero-point $x^{*}$ over the interval $[-1,1]$,
  $$x^{*}=\frac{1}{96} \left(63-\sqrt[3]{49041-16 \sqrt{3891895}}-\frac{1121}{\sqrt[3]{49041-16 \sqrt{3891895}}}\right)
  \approx -0.064041.$$
  Then $$(Z_5)_{\min}=\min\{Z_5(-1),Z_5(x^{*}),Z_5(1)\}=Z_5(x^{*})\approx14.3305.$$
  Thus, we have
  $$\mathrm{g}^{(5)}(\varphi)\ge \frac1{30}Z_5(x^{*})\approx0.477683.$$
  The proof is completed.
\end{proof}

\subsection{Decay estimates of DOC-$\rmk$ kernels}
\label{section: DOC decaying estimates}
\setcounter{equation}{0}

\begin{figure}[htb!]
\centering
\includegraphics[width=2in,height=1.6in]{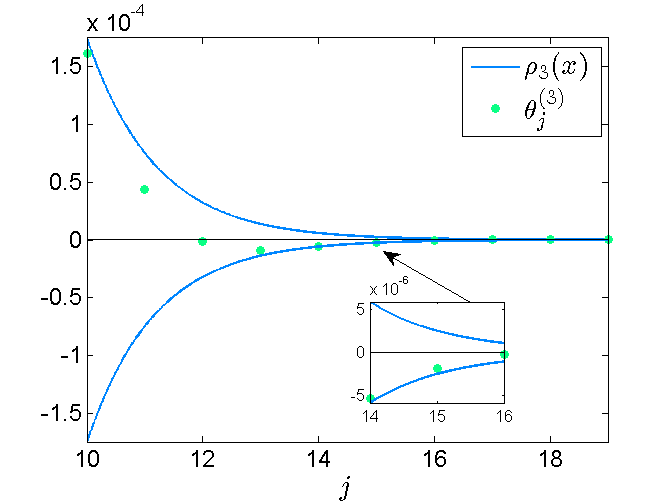}
\includegraphics[width=2in,height=1.6in]{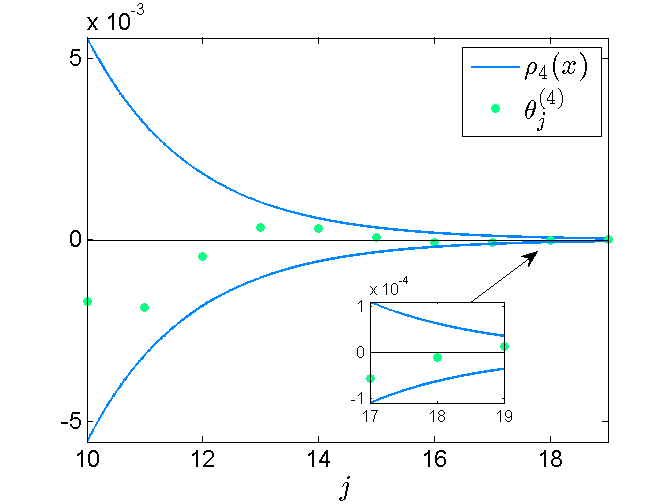}
\includegraphics[width=2in,height=1.6in]{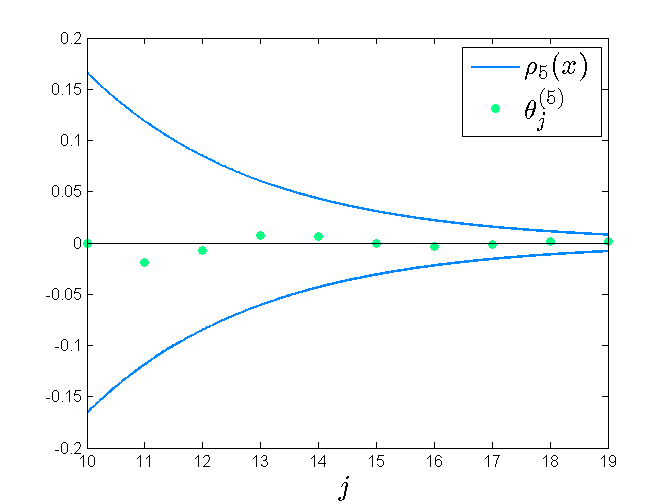}
\caption{The DOC-$\rmk$ kernels and the bound
$\rho_{\rmk}(x)=\frac{\rho_{\rmk}}{4}\bra{\frac{\rmk}{7}}^x$ for $3\le \rmk \le 5$.}
\label{fig: decaying property of DOC-k kernels}
\end{figure}

We now present the decay estimates of DOC-$\rmk$ kernels. Notice that although the BDF-$\rmk$ kernels $b_{j}^{(\rmk)}$ vanish for $j\ge \rmk$, the associated DOC-$\rmk$ kernels $\theta_{j}^{(\rmk)}$ are always nonzero for any $j\ge0$.
The following lemma presents the decay property of the DOC-$\rmk$ kernels (we plot in Figure \ref{fig: decaying property of DOC-k kernels} the decay properties of those kernels).

\begin{lemma}\label{lem: decaying estimates DOC}
The associated DOC-$\rmk$ kernels $\theta_{j}^{(\rmk)}$ defined in \eqref{def: DOC-Kernels} satisfy
\begin{align*}
\absb{\theta_{j}^{(\rmk)}}\le \frac{\rho_{\rmk}}{4}\braB{\frac{\rmk}{7}}^j\quad\text{for $3\le \rmk\le 5$ and $j\ge0$,}
\end{align*}
where $\rho_3=10/3$, $\rho_4=6$ and $\rho_5=96/5$.
\end{lemma}
\begin{proof}
By the definition \eqref{def: DOC-Kernels},
we have $\theta_{0}^{(\rmk)}=1/b_{0}^{(\rmk)}$ and
$\sum_{\ell=j}^{n}\theta_{n-\ell}^{(\rmk)}b_{\ell-j}^{(\rmk)}=0\,\, (\rmk\leq j \leq n-1),$
or,
\begin{align}\label{eq: difference equation DOC}
\sum_{m=0}^{j}\theta_{j-m}^{(\rmk)}b_{m}^{(\rmk)}=0
\quad\text{for $3\le \rmk\le 5$ and $j\ge0$.}
\end{align}
We will solve the difference equation \eqref{eq: difference equation DOC} to
find the solution $\theta_{j}^{(\rmk)}$ for any $j\ge0$.

\textbf{(I) The case $\rmk=3.$}\quad Taking $j=0$ and $j=1$ in \eqref{eq: difference equation DOC}
yield the two initial values
$$\theta_{0}^{(3)}=\frac{6}{11}\quad\text{and}\quad\theta_{1}^{(3)}=-\theta_{0}^{(3)}b_{1}^{(3)}/b_{0}^{(3)}=\frac{42}{11^2}.$$
One takes $j\ge2$ in \eqref{eq: difference equation DOC} and finds the equation
\begin{align*}
\theta_{j}^{(3)}b_{0}^{(3)}+\theta_{j-1}^{(3)}b_{1}^{(3)}+\theta_{j-2}^{(3)}b_{2}^{(3)}=0\quad\text{for $j\ge2$,}
\end{align*}
or
\begin{align}\label{eq: DOC-3 difference equation}
11\theta_{j}^{(3)}-7\theta_{j-1}^{(3)}+2\theta_{j-2}^{(3)}=0\quad\text{for $j\ge2$.}
\end{align}
The characteristic equation $11\lambda_3^2-7\lambda_3+2=0$ has two roots
$$\lambda_{3,1}=\bar{\lambda}_{3,2}=\frac{7+\imath\sqrt{39}}{22}.$$
Then it is easy to obtain the solution
\begin{align}\label{eq: DOC-3 formula}
\theta_{j}^{(3)}=\frac{39-7\imath\sqrt{39}}{143}\brab{\frac{7+\imath\sqrt{39}}{22}}^j
+\frac{39+7\imath\sqrt{39}}{143}\brab{\frac{7-\imath\sqrt{39}}{22}}^j\quad\text{for $j\ge0$,}
\end{align}
and the following decaying estimate
\begin{align}\label{eq: estimate DOC-3 formula}
\absb{\theta_{j}^{(3)}}\le 4\sqrt{\frac{6}{143}}\braB{\sqrt{\frac{2}{11}}}^{j}
\le \frac{5}{6}\braB{\frac{3}{7}}^{j}\quad\text{for $j\ge0$.}
\end{align}

\textbf{(II)  The case $\rmk=4$}\quad Taking $j=0,1$ and $j=2$ in \eqref{eq: difference equation DOC}
yield the initial values
$$\theta_{0}^{(4)}=\frac{12}{25},\quad\theta_{1}^{(4)}=\frac{12}{25}\frac{23}{25}
\quad\text{and}\quad\theta_{2}^{(4)}
=\frac{23}{25}\theta_{1}^{(4)}-\frac{13}{25}\theta_{0}^{(4)}
=\frac{48\cdot51}{25^3}.$$
One takes $j\ge3$ in \eqref{eq: difference equation DOC} and finds
\begin{align}\label{eq: DOC-4 difference equation}
\theta_{j}^{(4)}b_{0}^{(4)}+\theta_{j-1}^{(4)}b_{1}^{(4)}
+\theta_{j-2}^{(4)}b_{2}^{(4)}+\theta_{j-3}^{(4)}b_{3}^{(4)}=0\quad\text{for $j\ge3$.}
\end{align}
The characteristic equation $25\lambda_4^3-23\lambda_4^2+13\lambda_4-3=0$ has three roots
\begin{align*}
\lambda_{4,1}=&\,\bar{\lambda}_{4,2}=\frac{23}{75}-\frac{\brat{1+\imath \sqrt{3}}\nu}{75\sqrt[3]{4}}
+\frac{223 \brat{1-\imath \sqrt{3}}}{75\sqrt[3]{2}\nu}\approx 0.269261 - 0.492 \imath,\\
\lambda_{4,3}=&\,\frac{23}{75}+\frac{\sqrt[3]{2}\nu}{75}
-\frac{\sqrt[3]{4}}{75}\frac{223}{\nu}\approx0.381478,
\end{align*}
where $\nu:=\sqrt[3]{1921+225 \sqrt{511}}$ and
$$\absb{\lambda_{4,1}}=\absb{\lambda_{4,2}}\approx 0.560862.$$
We have the solution
\begin{align}\label{eq: DOC-4 formula}
\theta_{j}^{(4)}=d_{4,1}\lambda_{4,1}^j+d_{4,2}\lambda_{4,2}^j+d_{4,3}\lambda_{4,3}^j\quad\text{for $j\ge0$,}
\end{align}
where $d_{4,i}$ for $i=1,2,3$ are constant determined by the following equations
\begin{align*}
d_{4,1}\lambda_{4,1}^j+d_{4,2}\lambda_{4,2}^j+d_{4,3}\lambda_{4,3}^j=\theta_{j}^{(4)}\quad \text{for $j=0,1,2$.}
\end{align*}
Numerical computations yield
\begin{align*}
d_{4,1}=\bar{d}_{4,2}\approx -0.359522+0.405803\imath\quad\text{and}\quad  d_{4,3}\approx 0.719044.
\end{align*}
Then we can obtain the following estimate
\begin{align}\label{eq: estimate DOC-4 formula}
\absb{\theta_{j}^{(4)}}
\le\frac{3}{2}\braB{\frac{4}{7}}^j\quad\text{for $j\ge0$.}
\end{align}

\textbf{(III)  The case $\rmk=5$}\quad Taking $j=0,1,2$ and $j=3$ in \eqref{eq: difference equation DOC}
yield the initial values
\begin{align*}
\theta_{0}^{(5)}=&\,\frac{60}{137},\qquad\theta_{1}^{(5)}=\frac{163}{137}\theta_{0}^{(5)}=\frac{60\cdot163}{137^2},\\
\theta_{2}^{(5)}=&\,\frac{163}{137}\theta_{1}^{(5)}-\theta_{0}^{(5)}=\frac{60(163^2-1)}{137^3},\\
\theta_{3}^{(5)}=&\,\frac{163}{137}\theta_{2}^{(5)}-\theta_{1}^{(5)}+\frac{63}{137}\theta_{0}^{(5)}
=\frac{9780 (163^2-1)-6000\cdot137^2}{137^4}.
\end{align*}
One takes $j\ge4$ in \eqref{eq: difference equation DOC} and finds
\begin{align}\label{eq: DOC-5 difference equation}
\theta_{j}^{(5)}b_{0}^{(5)}+\theta_{j-1}^{(5)}b_{1}^{(5)}
+\theta_{j-2}^{(5)}b_{2}^{(5)}+\theta_{j-3}^{(5)}b_{3}^{(5)}
+\theta_{j-4}^{(5)}b_{4}^{(5)}=0\quad\text{for $j\ge4$.}
\end{align}
The characteristic equation $137\lambda_5^4-163\lambda_5^3+137\lambda_5^2-63\lambda_5+12=0$ has four roots
\begin{align*}
\lambda_{5,1}=\bar{\lambda}_{5,2}\approx0.210044 - 0.67687 \imath, \quad
\lambda_{5,3}=\bar{\lambda}_{5,4}\approx0.384847 - 0.162121 \imath.
\end{align*}
Also,
\begin{align*}
\absb{\lambda_{5,1}}=\absb{\lambda_{5,2}}\approx 0.708711 \quad\text{and}\quad
\absb{\lambda_{5,3}}=\absb{\lambda_{5,4}}\approx 0.417601
\end{align*}
We have the solution
\begin{align}\label{eq: DOC-5 formula}
\theta_{j}^{(5)}=\sum_{i=1}^4d_{5,i}\lambda_{5,i}^j\quad\text{for $j\ge0$.}
\end{align}
The constants $d_{5,i}$ for $i=1,2,3,4$ are determined by the following equations
\begin{align*}
\sum_{i=1}^4d_{5,i}\lambda_{5,i}^j=\theta_{j}^{(5)}\quad \text{for $j=0,1,2,3$,}
\end{align*}
which yield
\begin{align*}
d_{5,1}=\bar{d}_{5,2}\approx 0.0365741 - 0.450763 \imath, \quad d_{5,3}=\bar{d}_{5,4}\approx -0.0365741 + 3.27211 \imath.
\end{align*}
By using the fact
\begin{align*}
\abs{d_{5,1}}=\abs{d_{5,2}}\approx0.452244\quad\text{and}\quad \abs{d_{5,3}}=\abs{d_{5,4}}\approx3.27232,
\end{align*}
it is not difficult to obtain the following estimate
\begin{align}\label{eq: estimate DOC-5 formula}
\absb{\theta_{j}^{(5)}}
\le&\, \frac{24}{5}\braB{\frac{5}{7}}^j\quad\text{for $j\ge0$.}
\end{align}
This completes the proof.
\end{proof}

We close this section by noticing that the techniques in this section can also be used to handle the case $\rmk=2$ for which standard energy analysis is also applicable, see also in \cite{LiaoZhang:2020} for the analysis of variable time stepping.

\subsection{Decay Estimates of the starting values}
\label{section: Estimates of initial terms UI}
\setcounter{equation}{0}

Noticing that the starting values $u_{\mathrm{I}}^{(\rmk,j)}$
have different expressions \eqref{Dis: initial effect -BDF formula} for different step indexes $\rmk.$ We shall present the decay Estimates of the starting values by using
Lemma \ref{lem: decaying estimates DOC}.

\begin{lemma}\label{lem: Estimates of initial terms uI}
There exist positive constants $c_{\mathrm{I},\rmk}>1$ such that
the starting values $u_{\mathrm{I}}^{(\rmk,j)}$ satisfy
\begin{align*}
\absb{u_{\mathrm{I}}^{(\rmk,j)}}\le \frac{c_{\mathrm{I},\rmk}\rho_{\rmk}}{8}
\braB{\frac{\rmk}{7}}^{j-\rmk}\sum_{\ell=1}^{\rmk-1}
\absb{\diff u^\ell}\quad\text{for $3\le \rmk\le 5$ and $j\ge\rmk$,}
\end{align*}
such that
\begin{align*}
\sum_{j=\rmk}^n\absb{u_{\mathrm{I}}^{(\rmk,j)}}\le
\frac{7c_{\mathrm{I},\rmk}\rho_{\rmk}}{8(7-\rmk)}
\sum_{\ell=1}^{\rmk-1}
\absb{\diff u^\ell}\quad\text{for $3\le \rmk\le 5$ and $n\ge\rmk$,}
\end{align*}
where the constants $\rho_{\rmk}$ are defined in Lemma \ref{lem: decaying estimates DOC}.
\end{lemma}
\begin{proof}
\textbf{(I)  The case $\rmk=3$.}\quad
  Recalling the fact $b_{j}^{(3)}=0$ for $j\ge 3$, one can derive that
  \begin{align*}
 u_{\mathrm{I}}^{(3,n)}
=&\,\left\{
 \begin{array}{ll}
\theta_{0}^{(3)}b_{2}^{(3)}\diff u^1+\theta_{0}^{(3)}b_{1}^{(3)}\diff u^2, &\hbox{ for $n=3$,}\\[2ex]
\theta_{n-3}^{(3)}b_{2}^{(3)}\diff u^1+\brab{\theta_{n-3}^{(3)}b_{1}^{(3)}
+\theta_{n-4}^{(3)}b_{2}^{(3)}}\diff u^2, &\hbox{ for $n\ge4$}
\end{array}
\right.\\
=&\,\left\{
 \begin{array}{ll}
\frac{1}{11}\diff u^1-\frac{7}{11}\diff u^2, &\hbox{ for $n=3$,}\\[2ex]
\frac{1}{3}\theta_{n-3}^{(3)}\diff u^1-\frac{11}{6}\theta_{n-2}^{(3)}\diff u^2, &\hbox{ for $n\ge4$}
\end{array}
\right.
\end{align*}
 where the difference equation \eqref{eq: DOC-3 difference equation}
was used in the case of $n\ge4$.
So Lemma \ref{lem: decaying estimates DOC} yields
\begin{align*}
\absb{u_{\mathrm{I}}^{(3,3)}}\le&\,\frac{1}{11}\absb{\diff u^1}+\frac{7}{11}\absb{\diff u^2},\\
 \absb{u_{\mathrm{I}}^{(3,n)}}\le&\,\frac{1}{3}\absb{\theta_{n-3}^{(\rmk)}}\absb{\diff u^1}
 +\frac{11}{6}\absb{\theta_{n-2}^{(\rmk)}}\absb{\diff u^2}\\
 \le&\,\frac{\rho_3}{4}\braB{\frac{3}{7}}^{n-3}
\braB{\frac{1}{3}\absb{\diff u^1}+\frac{11}{14}\absb{\diff u^2}}\quad\text{for $n\ge4$.}
\end{align*}
The case $\rmk=3$ is verified by taking $c_{\mathrm{I},3}=11/7$ since
\begin{align*}
 \absb{u_{\mathrm{I}}^{(3,j)}}
\le&\,\frac{\rho_3}{4}\braB{\frac{3}{7}}^{j-3}
\braB{\frac{1}{3}\absb{\diff u^1}+\frac{11}{14}\absb{\diff u^2}}\quad\text{for $j\ge3$.}
\end{align*}

\textbf{(II)   The case $\rmk=4$.}\quad
It only needs to consider the general case $n\ge6$.
Since $b_{j}^{(4)}=0$ for $j\ge 4$, one has
\begin{align*}
 u_{\mathrm{I}}^{(4,n)}&\,=
\theta_{n-4}^{(4)}b_{3}^{(4)}\diff u^1
+\diff u^2\sum_{j=4}^{n}\theta_{n-j}^{(4)}b_{j-2}^{(4)}
+\diff u^3\sum_{j=4}^{n}\theta_{n-j}^{(4)}b_{j-3}^{(4)}\\
=&\,\theta_{n-4}^{(4)}b_{3}^{(4)}\diff u^1
+\brab{\theta_{n-4}^{(4)}b_{2}^{(4)}+\theta_{n-5}^{(4)}b_{3}^{(4)}}\diff u^2
-\theta_{n-3}^{(4)}b_{0}^{(4)}\diff u^3\quad\text{for  $n\ge6$,}
\end{align*}
where the difference equation \eqref{eq: DOC-4 difference equation} was used.
So Lemma \ref{lem: decaying estimates DOC} yields
\begin{align*}
\absb{u_{\mathrm{I}}^{(4,n)}}\le&\,
\frac{\rho_4}{4}\braB{\frac{4}{7}}^{n-4}\braB{\frac1{4}\absb{\diff u^1}
+\frac{3}{2}\absb{\diff u^2}
+\frac{25}{21}\absb{\diff u^3}}\quad\text{for $n\ge6$.}
\end{align*}
Then the estimate for $\rmk=4$ is verified by taking the fixed cases $n=4$ and 5 into account.

\textbf{(III)   The case $\rmk=5$.}\quad
We only consider the general case $n\ge8$. Since $b_{j}^{(5)}=0$ for $j\ge 5$, one has
\begin{align*}
 u_{\mathrm{I}}^{(5,n)}=\sum_{\ell=1}^{4}\diff u^\ell
\sum_{j=5}^{n}\theta_{n-j}^{(5)}b_{j-\ell}^{(5)}
=\theta_{n-5}^{(5)}b_{4}^{(5)}\diff u^1
+\brab{\theta_{n-5}^{(5)}b_{3}^{(5)}+\theta_{n-6}^{(5)}b_{4}^{(5)}}\diff u^2\\
+\brab{\theta_{n-5}^{(5)}b_{2}^{(5)}+
\theta_{n-6}^{(5)}b_{3}^{(5)}+\theta_{n-7}^{(5)}b_{4}^{(5)}}\diff u^3
-\theta_{n-4}^{(5)}b_{0}^{(5)}\diff u^4\quad\text{for $n\ge8$,}
\end{align*}
where the difference equation \eqref{eq: DOC-5 difference equation} was used
in the last term. By using Lemma \ref{lem: decaying estimates DOC}, one gets
\begin{align*}
\absb{u_{\mathrm{I}}^{(5,n)}}\le&\,
\frac{\rho_5}{4}\braB{\frac{5}{7}}^{n-5}\braB{\frac1{5}\absb{\diff u^1}
+\frac{4}{3}\absb{\diff u^2}
+\frac{25}{6}\absb{\diff u^3}+\frac{23}{14}\absb{\diff u^4}}\quad\text{for $n\ge8$.}
\end{align*}
The estimate for $\rmk=5$ can be  verified with a finite $c_{\mathrm{I},5}$
by taking the fixed cases $n=5, 6$ and 7 into account. The proof is completed.
\end{proof}

\section{Discrete energy analysis for linear reaction-diffusion}

We are now ready to present the main results of this work.

\subsection{Stability analysis}
\setcounter{equation}{0}
We first consider the dissipative case with $\beta=\beta(x)\le0.$ In this case, we have the following stability result.

\begin{theorem}\label{thm: L2 norm stability}
The time-discrete solution $u^n$ of the BDF-$\rmk$ $(3\le \rmk\le 5)$ scheme
\eqref{eq: time-discrete IBVP}
for the dissipative case $\beta=\beta(x)\le0$ satisfies
\begin{align*}
\mynormb{u^{n}}
\le&\,\mynormb{u^{\rmk-1}}+\frac{7c_{\mathrm{I},\rmk}\rho_{\rmk}}{4(7-\rmk)}\sum_{\ell=1}^{\rmk-1}\mynormb{\diff u^\ell}
+\frac{7\rho_{\rmk}}{2(7-\rmk)}\sum_{\ell=\rmk}^n\tau\mynormb{f^\ell}\\
\le&\,\frac{7\rho_{\rmk}}{2(7-\rmk)}\braB{c_{\mathrm{I},\rmk}\sum_{\ell=0}^{\rmk-1}\mynormb{u^\ell}
+\sum_{\ell=\rmk}^n\tau\mynormb{f^\ell}}
\qquad\text{for $n\ge\rmk$,}
\end{align*}
where the constants $\rho_{\rmk}$ and $c_{\mathrm{I},\rmk}$ are defined in Lemmas \ref{lem: decaying estimates DOC}
and \ref{lem: Estimates of initial terms uI}, respectively.
\end{theorem}
\begin{proof}Lemmas \ref{lem: equvielant positive definite DOC}
and \ref{lem: BDF-k minimum eigenvalue}
in Appendix \ref{section: DOC positive definiteness} imply that
the DOC-$\rmk$ kernels are positive definite for $3\le \rmk\le 5$. Under the setting $\beta=\beta(x),$ one has
\begin{align*}
&-2\varepsilon\sum_{j=\rmk}^n\sum_{\ell=\rmk}^j
\theta_{j-\ell}^{(\rmk)}\myinnerb{\nabla u^{\ell},\nabla u^j}\le0\quad\text{and}\quad
2\sum_{j=\rmk}^n\sum_{\ell=\rmk}^j\theta_{j-\ell}^{(\rmk)}\myinnerb{\beta  u^{\ell},u^j}\le0.
\end{align*}
It follows from \eqref{eq: DOC IBVP energy form} that
\begin{align}\label{eq: DOC IBVP energy form-1A}
\mynormb{u^{n}}^2\le&\,\mynormb{u^{\rmk-1}}^2-2\sum_{j=\rmk}^n\myinnerb{u_{\mathrm{I}}^{(\rmk,j)},u^j}
        +2\tau\sum_{j=\rmk}^n\sum_{\ell=\rmk}^j\myinnerb{\theta_{j-\ell}^{(\rmk)}f^\ell,u^j}\nonumber\\
        \le&\,\mynormb{u^{\rmk-1}}^2+2\sum_{j=\rmk}^n\mynormb{u_{\mathrm{I}}^{(\rmk,j)}}\mynormb{u^j}
        +2\tau\sum_{j=\rmk}^n\sum_{\ell=\rmk}^j\mynormb{\theta_{j-\ell}^{(\rmk)}f^\ell}\mynormb{u^j}\quad\text{for $n\ge\rmk$}.
\end{align}
Taking some integer $n_1$ ($\rmk-1\le n_1\le n$) such that $\mynormb{u^{n_1}}=\max_{\rmk-1\le j\le n}\mynormb{u^{j}}$. Taking $n:=n_1$ in the above inequality, one gets
\begin{align*}
\mynormb{u^{n_1}}^2\le  \mynormb{u^{\rmk-1}}\mynormb{u^{n_1}}
+2\mynormb{u^{n_1}}\sum_{j=\rmk}^{n_1}\mynormb{u_{\mathrm{I}}^{(\rmk,j)}}
+2\tau\mynormb{u^{n_1}}\sum_{j=\rmk}^{n_1}\sum_{\ell=\rmk}^j\mynormb{\theta_{j-\ell}^{(\rmk)}f^\ell},
\end{align*}
and thus
\begin{align}\label{eq: DOC IBVP energy form-1}
\mynormb{u^{n}}\le&\,\mynormb{u^{n_1}}\le  \mynormb{u^{\rmk-1}}+2\sum_{j=\rmk}^{n_1}\mynormb{u_{\mathrm{I}}^{(\rmk,j)}}
+2\tau\sum_{j=\rmk}^{n_1}\sum_{\ell=\rmk}^j\mynormb{\theta_{j-\ell}^{(\rmk)}f^\ell}\nonumber\\
\le&\, \mynormb{u^{\rmk-1}}+2\sum_{j=\rmk}^{n}\mynormb{u_{\mathrm{I}}^{(\rmk,j)}}
+2\tau\sum_{j=\rmk}^{n}\sum_{\ell=\rmk}^j\mynormb{\theta_{j-\ell}^{(\rmk)}f^\ell}\quad\text{for $n\ge\rmk$.}
\end{align}
Applying Lemma \ref{lem: decaying estimates DOC},
we have $\absb{\theta_{j-\ell}^{(\rmk)}}\le \frac{\rho_{\rmk}}{4}(\frac{\rmk}{7})^{j-\ell}$ for $3\le \rmk\le 5$ and then
\begin{align}\label{ieq: exterior force estimate}
2\tau\sum_{j=\rmk}^{n}\sum_{\ell=\rmk}^j\mynormb{\theta_{j-\ell}^{(\rmk)}f^\ell}
\le&\, 2\tau\sum_{j=\rmk}^{n}\sum_{\ell=\rmk}^j\absb{\theta_{j-\ell}^{(\rmk)}}\mynormb{f^\ell}
=2\tau\sum_{\ell=\rmk}^n\mynormb{f^\ell}\sum_{j=\ell}^{n}\absb{\theta_{j-\ell}^{(\rmk)}}\nonumber\\
\leq&\, \frac{\rho_{\rmk}}{2}\sum_{\ell=\rmk}^n\tau\mynormb{f^\ell}\sum_{j=\ell}^{n}(\rmk/7)^{j-\ell}
\le \frac{7\rho_{\rmk}}{2(7-\rmk)}\sum_{\ell=\rmk}^n\tau\mynormb{f^\ell}
\quad\text{for $n\ge\rmk$.}
\end{align}
Then the claimed estimate follows by using together \eqref{eq: DOC IBVP energy form-1} and Lemma \ref{lem: Estimates of initial terms uI}.
\end{proof}

Next we consider the general case $\abs{\beta(x,t)}\le \beta^{*}$. In this case, we have following stability result.

\begin{theorem}\label{thm: general L2 norm stability}
Consider $3\le \rmk\le 5$ and the bounded coefficient $\beta(x,t)\le \beta^{*}$.
If the time-step size $\tau\le(7-\rmk)/(7\rho_{\rmk}\beta^{*})$,
the time-discrete solution $u^n$ of the BDF-$\rmk$
scheme \eqref{eq: time-discrete IBVP} satisfies
\begin{align*}
\mynormb{u^{n}}\le&\, \frac{7\rho_{\rmk}}{7-\rmk}\exp\braB{\frac{7\rho_{\rmk}}{7-\rmk}\beta^{*}t_{n-\rmk}}
\braB{c_{\mathrm{I},\rmk}\sum_{\ell=0}^{\rmk-1}\mynormb{u^\ell}
+\sum_{\ell=\rmk}^n\tau\mynormb{f^\ell}}\quad\text{for $\rmk\le n\le N$}.
\end{align*}
where the constants $\rho_{\rmk}$ and $c_{\mathrm{I},\rmk}$ are defined by Lemmas \ref{lem: decaying estimates DOC}
and \ref{lem: Estimates of initial terms uI}, respectively.
\end{theorem}
\begin{proof}
By the inequality \eqref{eq: DOC IBVP energy form}
and Theorem \ref{thm: L2 norm stability}, we aim at bounding the third term
of the right hand side of \eqref{eq: DOC IBVP energy form},
\begin{align}\label{ieq: direct estimate reaction term}
2\tau\sum_{j=\rmk}^n\sum_{\ell=\rmk}^j\theta_{j-\ell}^{(\rmk)}\myinnerb{\beta^\ell u^{\ell},u^j}
\le&\,2\beta^{*}\tau\sum_{j=\rmk}^n\sum_{\ell=\rmk}^j\absb{\theta_{j-\ell}^{(\rmk)}}\mynormb{u^\ell}\mynormb{u^j}.
\end{align}
Then it is not difficult to derive from \eqref{eq: DOC IBVP energy form} that,
for $\rmk\le n\le N$,
\begin{align}\label{ieq: DOC IBVP energy form-beta1}
\mynormb{u^{n}}^2
        \le\mynormb{u^{\rmk-1}}^2&\,+2\sum_{j=\rmk}^n\mynormb{u_{\mathrm{I}}^{(\rmk,j)}}\mynormb{u^j}
+2\tau\sum_{j=\rmk}^n\sum_{\ell=\rmk}^j
\braB{\beta^{*}\absb{\theta_{j-\ell}^{(\rmk)}}\mynormb{u^\ell}+\mynormb{\theta_{j-\ell}^{(\rmk)}f^\ell}}\mynormb{u^j}.
\end{align}
Taking some integer $n_2$ ($\rmk-1\le n_2\le n$) such that $\mynormb{u^{n_2}}=\max_{\rmk-1\le j\le n}\mynormb{u^{j}},$ and setting $n:=n_2$ in the above inequality \eqref{ieq: DOC IBVP energy form-beta1}, one obtains
\begin{align*}
\mynormb{u^{n_2}}\le  \mynormb{u^{\rmk-1}}
+2\sum_{j=\rmk}^{n_2}\mynormb{u_{\mathrm{I}}^{(\rmk,j)}}
+2\beta^{*}\tau\sum_{j=\rmk}^{n_2}\mynormb{u^j}\sum_{\ell=\rmk}^j\absb{\theta_{j-\ell}^{(\rmk)}}
+2\tau\sum_{j=\rmk}^{n_2}\sum_{\ell=\rmk}^j\mynormb{\theta_{j-\ell}^{(\rmk)}f^\ell},
\end{align*}
and thus
\begin{align*}
\mynormb{u^{n}}\le  \mynormb{u^{\rmk-1}}
+2\sum_{j=\rmk}^{n}\mynormb{u_{\mathrm{I}}^{(\rmk,j)}}
+2\beta^{*}\tau\sum_{j=\rmk}^{n}\mynormb{u^j}\sum_{\ell=\rmk}^j\absb{\theta_{j-\ell}^{(\rmk)}}
+2\tau\sum_{j=\rmk}^{n}\sum_{\ell=\rmk}^j\mynormb{\theta_{j-\ell}^{(\rmk)}f^\ell},
\end{align*}
By applying Lemma \ref{lem: decaying estimates DOC} we have
$$\sum_{\ell=\rmk}^j\absb{\theta_{j-\ell}^{(\rmk)}}
\le \frac{\rho_{\rmk}}{4}\sum_{\ell=\rmk}^j(\rmk/7)^{j-\ell}
\le \frac{7\rho_{\rmk}}{4(7-\rmk)}\quad \text{for $3\le \rmk\le 5$.}$$
Then we apply Lemma \ref{lem: Estimates of initial terms uI} and the estimate \eqref{ieq: exterior force estimate} to find that
\begin{align}\label{ieq: DOC IBVP energy form-beta2}
\mynormb{u^{n}}\le&\,  \mynormb{u^{\rmk-1}}
+\frac{7c_{\mathrm{I},\rmk}\rho_{\rmk}}{4(7-\rmk)}\sum_{\ell=1}^{\rmk-1}\mynormb{\diff u^\ell}
+\frac{7\rho_{\rmk}\beta^{*}}{2(7-\rmk)}\sum_{j=\rmk}^{n}\tau\mynormb{u^j}
+\frac{7\rho_{\rmk}}{2(7-\rmk)}\sum_{\ell=\rmk}^n\tau\mynormb{f^\ell}\nonumber\\
\le&\, \frac{7c_{\mathrm{I},\rmk}\rho_{\rmk}}{2(7-\rmk)}\sum_{\ell=0}^{\rmk-1}\mynormb{u^\ell}
+\frac{7\rho_{\rmk}\beta^{*}}{2(7-\rmk)}\sum_{j=\rmk}^{n}\tau\mynormb{u^j}
+\frac{7\rho_{\rmk}}{2(7-\rmk)}\sum_{\ell=\rmk}^n\tau\mynormb{f^\ell}.
\end{align}
If the time-step size $\tau\le\frac{7-\rmk}{7\rho_{\rmk}\beta^{*}}$,
it follows from \eqref{ieq: DOC IBVP energy form-beta2} that
\begin{align*}
\mynormb{u^{n}}\le&\, \frac{7c_{\mathrm{I},\rmk}\rho_{\rmk}}{7-\rmk}\sum_{\ell=0}^{\rmk-1}\mynormb{u^\ell}
+\frac{7\rho_{\rmk}\beta^{*}}{7-\rmk}\sum_{j=\rmk}^{n-1}\tau\mynormb{u^j}
+\frac{7\rho_{\rmk}}{7-\rmk}\sum_{\ell=\rmk}^n\tau\mynormb{f^\ell}\quad\text{for $\rmk\le n\le N$}.
\end{align*}
Then the claimed estimate follows by using the standard Gr\"{o}nwall inequality.
\end{proof}

\subsection{Convergence analysis}

Let $\tilde{u}^n:=u(t_n,x)-u^n(x)$ for $n\ge0$.
Then the error equation of \eqref{eq: time-discrete IBVP} reads
\begin{align}\label{eq: time-discrete IBVP error}
D_{\rmk}\tilde{u}^{n}
	=\varepsilon\Delta \tilde{u}^{n}
        +\beta^n \tilde{u}^{n}+\eta^n,\quad\text{for $\rmk\le n\le N$,}
\end{align}
where the local consistency error $\eta^j=D_{\rmk}u(t_j)-\partial_tu(t_j)$ for $j\ge\rmk$.
Assume that the solution is regular in time for $t\ge t_{\rmk}$ such that
$$\absb{\eta^j}\le C_u\tau^{\rmk}\max_{t_{\rmk}\le t\le T}\absb{\partial_t^{(\rmk+1)}u(t)}
\le C_u\tau^{\rmk}\quad
\text{for $j\ge\rmk$}.$$
The stability estimate
in Theorem \ref{thm: general L2 norm stability}  yields
\begin{align*}
\mynormb{\tilde{u}^{n}}\le&\, \frac{7\rho_{\rmk}}{7-\rmk}
\exp\braB{\frac{7\rho_{\rmk}}{7-\rmk}\beta^{*}t_{n-\rmk}}
\braB{c_{\mathrm{I},\rmk}\sum_{\ell=0}^{\rmk-1}\mynormb{\tilde{u}^\ell}
+\sum_{\ell=\rmk}^n\tau\mynormb{\eta^\ell}}\quad\text{for $\rmk\le n\le N$}.
\end{align*}
This implies at the following theorem.
\begin{theorem}\label{thm: general L2 norm converegnce}
Let $u(t_n,x)$ and $u^n(x)$ be the solutions of
the diffusion problem \eqref{eq: diffusion problem}
and the BDF-$\rmk$ scheme \eqref{eq: time-discrete IBVP}, respectively.
If the time-step size $\tau\le(7-\rmk)/(7\rho_{\rmk}\beta^{*})$,
then the time-discrete solution $u^n$ is convergent in the $L^2$ norm,
\begin{align*}
\mynormb{u(t_n)-u^{n}}\le&\,
\frac{7\rho_{\rmk}c_{\mathrm{I},\rmk}}{7-\rmk}
\exp\braB{\frac{7\rho_{\rmk}\beta^{*}t_{n-\rmk}}{7-\rmk}}
\bigg(\sum_{\ell=0}^{\rmk-1}\mynormb{u(t_{\ell})-u^{\ell}}
+C_ut_{n-\rmk}\tau^{\rmk}\bigg)
\end{align*}
for $\rmk\le n\le N$, where
$\rho_{\rmk}$ and $c_{\mathrm{I},\rmk}$ are defined by Lemmas \ref{lem: decaying estimates DOC}
and \ref{lem: Estimates of initial terms uI}, respectively.
\end{theorem}

\section{Concluding remarks}

In this work, we presented a novel discrete energy analysis for the BDF-$\rmk$ schemes with $3\le \rmk\le 5$ by using the \textit{discrete orthogonal convolution kernels} technique. Our analysis is straightforward in the sense that the standard inner product with $u^j$ is adopted, which coincides with the classical energy approach. With this straightforward approach, less spacial regularity requirement is required for the initial data by comparing with the multiplier technique which requires stronger norm for the initial data, see the stability estimates in \cite[Proposition 5.1 and Theorem 5.1]{AkrivisKatsoprinakis:2016}.

There are several remaining issues to be handled in future works:
\begin{itemize}
\item
 The present work opens up the possibility for handling the BDF-$\rmk$ time-stepping methods for nonlinear diffusion problems, such as phase filed models \cite{LiaoTangZhou:2020bdf2}.
 \item
  Note that a very recent work by Akrivis et al. \cite{AkrivisChenYuZhou:2021} successfully applied the multiplier technique to deal with the BDF-6 scheme. It is known that a linear multistep method is {\em zero-stable} if a perturbation in the starting values of size $\epsilon$ causes the numerical solution over any time interval to change by no more than $O(\epsilon)$. This is called {\em zero-stability} because it is enough to check the condition for the differential equation $y' = 0$. It is known that the BDF-$\rmk$ formulas for $\rmk >6$ are not zero-stable so they cannot be used. Note that the current analysis can only cover $\rmk$ up to 5, and finer analysis is needed to extend the DOC kernels technique for the BDF-6 scheme.
   \item
  Another interesting topic is to investigate the discrete energy technique for the stability and convergence
  of the BDF-$\rmk$ formula ($3\le \rmk\le 5$) with variable time-steps.
\end{itemize}


\begin{thebibliography}{99}

\bibitem{Akrivis:2015}
{\sc G. Akrivis},
{\em Stability of implicit-explicit backward difference
formulas for nonlinear parabolic equations},
SIAM J. Numer. Anal., 53 (2015), pp. 464-484.



\bibitem{AkrivisChenYuZhou:2021}
{\sc G. Akrivis, M.H. Chen, F. Yu, Z. Zhou},
{\em The energy technique for the six-step BDF method},
Math. Comp., 2021, Revised, arXiv:2007.08924.

\bibitem{AkrivisKatsoprinakis:2016}
{\sc G. Akrivis, E. Katsoprinakis},
{\em Backward difference formulae: new multipliers and stability properties
for parabolic equations}, Math. Comp.,  85 (2016), pp. 2195-2216.


\bibitem{AkrivisLubich:2015}
{\sc G. Akrivis and C. Lubich},
{\em Fully implicit, linearly implicit and implicit-explicit backward difference
formulae for quasi-linear parabolic equations},
Numer. Math., 131 (2015), pp. 713-735.


\bibitem{Dahlquist:1978}
{\sc G. Dahlquist}, {\em G-stability is equivalent to A-stability},
BIT, 18 (1978), pp. 384-401.

\bibitem{GrenanderSzego:2001}
{\sc U. Grenander and G. Szeg\"{o}}, {\em Toeplitz Forms and Their Applications},
2nd edition, AMS Chelsea, Providence, RI, 2001.

\bibitem{HairerNorsettWanner:2002}
{\sc E. Hairer and G. Wanner},
{\em Solving Ordinary Differential Equations II:
Stiff and Differential-Algebraic Problems},
Springer Series in Computational Mathematics Volume 14,
Second Edition, Springer-Verlag, 2002.


%
%


\bibitem{LiaoTangZhou:2020bdf2}
{\sc H.-L. Liao, T. Tang and T. Zhou},
{\em On energy stable,  maximum-bound preserving, second-order BDF scheme
with variable steps for the Allen-Cahn equation},
SIAM J. Numer. Anal., 58(4) (2020), pp. 2294-2314.


\bibitem{LiaoTangZhou:2020doc}
{\sc H.-L. Liao, T. Tang and T. Zhou},
{\em Positive definiteness of real quadratic forms resulting
from variable-step approximations of convolution operators}, arXiv:2011.13383v1,
2020.

\bibitem{LiaoZhang:2020}
{\sc H.-L. Liao and Z. Zhang},
{\em Analysis of adaptive BDF2 scheme for diffusion equations},
Math. Comp., 2020, DOI: 10.1090/mcom/3585.


\bibitem{Liu:2013}
{\sc J. Liu},
{\em Simple and efficient ALE methods with provable
temporal accuracy up to fifth order for the
stokes equations on time varying domains}.
SIAM J. Numer. Anal., 51 (2013), pp. 743-772.

\bibitem{LubichMansourVenkataraman:2013}
{\sc C. Lubich, D. Mansour, and C. Venkataraman},
{\em Backward difference time discretization of parabolic
differential equations on evolving surfaces},
IMA J. Numer. Anal., 33 (2013), pp. 1365-1385.

\bibitem{NevanlinnaOdeh:1981}
{\sc O. Nevanlinna and F. Odeh},
{\em Multiplier techniques for linear multistep methods},
Numer. Funct. Anal. Optim. 3 (1981), pp. 377-423.


\bibitem{Thomee:2006}
{\sc V. Thom\'{e}e},
{\em Galerkin Finite Element Methods for Parabolic Problems},
Second Edition,  Springer-Verlag, 2006.


\end{thebibliography}
\end{document}